\documentclass[9pt,amsfonts]{article}
\usepackage{graphicx,amssymb,eucal,mathrsfs}
\usepackage{latexsym,amsmath,amscd,epsfig,xcolor}
\usepackage[all]{xy}
\usepackage{ntheorem,tikz}
\usepackage{epsf,graphicx,pgf,etex,amscd,pstricks}
\usepackage[french,english]{babel}
\usepackage[hypertexnames=false,backref=page,pdftex,
 	pdfpagemode=UseNone,
 	breaklinks=true,
 	extension=pdf,
 	colorlinks=true,
 	linkcolor=blue,
 	citecolor=red,
 	urlcolor=blue,
 ]{hyperref}

\newcommand{\C}{\mathbb{C}}

\newtheorem{main}{Theorem}

\newtheorem{prop}[main]{Proposition}

\newtheorem{rem}[main]{Remark}

\newtheorem{imain}{Theorem}

\theoremstyle{definition}

\title{Irrational pencils and Betti numbers}
\date{January 2021}
\author{Francisco Nicol\'as and Pierre Py$^{{\rm a}}$}

\begin{document}

\selectlanguage{english}

\maketitle

\begin{abstract} We study irrational pencils with isolated critical points on compact aspherical complex manifolds. We prove that if the set of critical points is nonempty, the homology of the kernel of the morphism induced by the pencil on fundamental groups is not finitely generated. This generalizes a result by Dimca, Papadima and Suciu. By considering self-products of the Cartwright-Steger surface, this allows us to build new examples of smooth projective varieties whose fundamental group has a non-finitely generated homology. 
\end{abstract}

\selectlanguage{french}

\begin{abstract} Nous \'etudions les pinceaux irrationnels \`a points critiques isol\'es sur les vari\'et\'es complexes compactes et asph\'eriques. Nous prouvons que si un tel pinceau poss\`ede au moins un point critique, alors l'homologie du noyau du morphisme induit entre groupes fondamentaux n'est pas de type fini. Ceci g\'en\'eralise un r\'esultat de Dimca, Papadima et Suciu. En consid\'erant le produit de plusieurs copies de la surface de Cartwright-Steger, ceci nous permet de donner de nouveaux exemples de vari\'et\'es projectives lisses dont le groupe fondamental a un groupe d'homologie qui n'est pas de type fini. 
\end{abstract}

\selectlanguage{english}

\footnotetext[1]{Partially supported by the french project ANR AGIRA.}

\section{Introduction}

\subsection{Finiteness properties and fundamental groups of smooth projective varieties}

Recall that a group $G$ is of type $\mathscr{F}_{\ell}$, for some integer $\ell$, if there exists a classifying space for $G$ (i.e. a $K(G,1)$) which is a ${\rm CW}$-complex with finite $\ell$-skeleton. This condition was introduced by Wall~\cite{wall}. For finitely presented groups, this is equivalent to the property of being of type ${\rm FP}_{\ell}$ (see \cite{brown} for the definition of this last property). Given a group $G$ one can ask whether it admits a classifying space which is a finite complex or whether it has property $\mathscr{F}_{\ell}$ for some $\ell$. We refer to~\cite{bb,bieri,bhms,stallings} for important works on these notions. 

In the context of the study of fundamental groups of smooth projective varieties, called {\it projective groups} in what follows, Koll\'ar asked in~\cite[\S 0.3.1]{kollar} whether a projective group is always commensurable (up to finite kernels) to a group admitting a classifying space which is a quasi-projective variety. Since any quasi-projective variety has the homotopy type of a finite complex~\cite[p. 27]{dimcalivre}, a positive answer to this question would imply that any projective group is commensurable to a group having a finite classifying space. However, a negative answer to Koll\'ar's question was given by Dimca, Papadima and Suciu in~\cite{dps}. To describe their results let us introduce some notations. 

{\bf Notations.} Throughout this text, $X$ will be a (connected) compact complex manifold of complex dimension $n\ge 2$ and $S$ will be a closed Riemann surface of positive genus. A surjective holomorphic map with connected fibers $f : X \to S$ is called an {\it irrational pencil}. For such a map, we will always denote by $\Lambda$ the kernel of the induced homomorphism $f_{\ast} : \pi_{1}(X)\to \pi_{1}(S)$.

In~\cite{dps} the authors proved the following two results. 

\begin{main}[see Theorem C in~\cite{dps}]\label{dpsfirst}
{\it If $n\ge 3$ and if $f : X\to S$ is an irrational pencil with isolated critical points, then the fundamental group of a smooth fiber of $f$ embeds into that of $X$ and coincides with the kernel $\Lambda$ of the induced homomorphism $f_{\ast} : \pi_{1}(X)\to \pi_{1}(S)$.}
\end{main}

\begin{main}[see \S 2 in~\cite{dps}]\label{dpsbetti}
{\it Let $X=\Sigma_{1}\times \cdots \times \Sigma_{n}$ be a direct product of $n$ Riemann surfaces of genus greater than $1$ and let $S$ have genus $1$. If $f : X \to S$ is an irrational pencil with isolated critical points then the group $H_{n}(\Lambda, \mathbb{Q})$ is infinite dimensional.} 
\end{main}

\begin{rem}
If we use the group structure on $S$, it is easy to see that any irrational pencil $f : X \to S$ as in Theorem \ref{dpsbetti} is the sum of holomorphic maps $f_i : \Sigma_{i}\to S$, i.e. $$f(p_1,\ldots , p_n)=f_{1}(p_{1})+\cdots +f_{n}(p_{n}).$$ 
If all the $f_i$'s are nonconstant and $n\ge 2$, $f$ has connected fibers if and only if it is $\pi_{1}$-surjective, see Lemma 2.1 in~\cite{Llo-fourier}. The set of critical points of $f$ is the product of the critical sets of the $f_i$'s. Hence $f$ has isolated critical points if and only if all the $f_i$'s are nonconstant.
\end{rem}

Combining Theorems \ref{dpsfirst} and \ref{dpsbetti}, Dimca, Papadima and Suciu answered negatively Koll\'ar's question. Indeed in the situation of Theorem~\ref{dpsbetti} the group $\Lambda$, which is the fundamental group of a smooth fiber of $f$ if $n\ge 3$, cannot be of type ${\rm FP}_{n}$ as the group $H_{n}(\Lambda,\mathbb{Z})$ is not finitely generated. The property of being of type ${\rm FP}_{n}$ is invariant by the commensurability relation~\cite{bieri,dps}. Therefore, no group commensurable to $\Lambda$ can have a finite classifying space.   

Building on the work~\cite{dps}, further examples of fundamental groups of smooth projective varieties with exotic finiteness properties were constructed and studied by Llosa Isenrich~\cite{Llo-fourier,Llo-17} and by Bridson and Llosa Isenrich~\cite{bli}. Note however that all the examples studied in~\cite{bli,Llo-fourier,Llo-17} are either subgroups of direct products of surface groups or extensions of such subgroups as in~\cite{bli}. The purpose of this note is to provide further examples of projective groups violating the $\mathscr{F}_{\ell}$ condition for some $\ell$ and which are not isomorphic to subgroups of direct products of surface groups (see Proposition~\ref{paspsg}), as well as to prove that the conclusion of Theorem~\ref{dpsbetti} holds in much greater generality. 

\subsection{Irrational pencils and Betti numbers}\label{ipabn}

Let $f : X \to S$ be an irrational pencil with ${\rm dim}_{\mathbb{C}} X=n \geq 2$. As before $\Lambda$ denotes the kernel of the morphism $f_{\ast} : \pi_{1}(X)\to \pi_{1}(S)$. We assume that the critical points of $f$ are isolated and that $f$ is not a submersion; its critical set is then a nonempty finite set. Let $\widehat{X} \to X$ be the covering space such that $\pi_{1}(\widehat{X})\simeq \Lambda$. Our main results are the following: 

\begin{imain}\label{bnestinfini}
The homology group $H_{n}(\widehat{X},\mathbb{Q})$ is infinite dimensional. 
\end{imain}

\begin{imain}\label{thgroup} If $X$ is aspherical, the group $H_{n}(\Lambda, \mathbb{Q})$ is infinite dimensional. In particular $\Lambda$ is not of type ${\rm FP}_{n}$. 
\end{imain}

Note that if $X$ is aspherical, the space $\widehat{X}$ is a $K(\Lambda,1)$, hence Theorem~\ref{thgroup} follows from Theorem~\ref{bnestinfini}. The very short proof of Theorem~\ref{bnestinfini} will be given in Section~\ref{proof}. 

When $n\ge 3$, the group $\Lambda$ is isomorphic to the fundamental group of a generic fiber of $f$, thanks to Theorem~\ref{dpsfirst}. Hence, when $n\ge 3$ and $X$ is aspherical and projective, $\Lambda$ is a projective group with exotic finiteness properties. We now discuss the relation of our results with the earlier works~\cite{bmp,dps,kapovich,Llo-fourier,Llo-17}.
 
\noindent $\bullet$ To study the topology of the covering space $\widehat{X}$, we follow the same method as in the articles~\cite{dps,kapovich}, which can be seen as a complex analog of Bestvina and Brady's work~\cite{bb}. Namely, denoting by $\widehat{S}$ the universal cover of $S$, we consider a lift $\widehat{f} : \widehat{X} \to \widehat{S}$ of $f$ and analyze the topology of $\widehat{X}$ by viewing it as the union $\cup_{k\ge 1}g^{-1}(D_{k})$ where $(D_{k})_{k\ge 1}$ is an increasing union of disks in $\widehat{S}$ and $g : \widehat{X}\to \widehat{S}$ is either the map $\widehat{f}$ or a small perturbation of it. 
 
\noindent $\bullet$ In the special case when $X$ is a product of Riemann surfaces and $S$ has genus $1$, Theorem~\ref{thgroup} reduces to Theorem~\ref{dpsbetti} above. Dimca, Papadima and Suciu's proof was based on the notion of {\it characteristic varieties}. Our proof however is more direct and applies in full generality: $X$ can be any aspherical complex manifold and the genus of $S$ need not be equal to $1$\footnote{On the other hand the article~\cite{dps} also studies the finiteness properties of arbitrary normal coabelian subgroups of direct products of fundamental groups of closed surfaces, not necessarily coming from irrational pencils.}. 

\noindent $\bullet$ When $X$ is a product of Riemann surfaces, the fact that $\Lambda$ is not of type ${\rm FP}_{n}$ can be deduced from the work of Bridson, Howie, Miller and Short~\cite{bhms}. See~\cite{Llo-fourier,Llo-17} for further results which rely on properties of subgroups of direct products of surface groups. 

\noindent $\bullet$ Under the same assumptions as in Theorem~\ref{thgroup}, the article~\cite{bmp} proves the weaker result that the group $\Lambda$ is not of type ${\rm FP}$ see~\cite[ \S VIII.6]{brown} for the definition of this property and Theorem 7.3 in~\cite{bmp} for this result.  

\noindent $\bullet$ When $n=2$, Kapovich~\cite{kapovich} has proved that if $f : X \to S$ is an irrational pencil which is not a submersion, with no multiple fiber and with $X$ aspherical, then the group $\Lambda=\ker(f_{\ast})$ is not finitely presented (without assuming that $f$ has isolated critical points). In the case of maps with isolated critical points, Theorem~\ref{thgroup} gives a slight strengthening of Kapovich's result since a finitely presented group must have finitely generated second homology group. Our proof is very close in spirit to the one in~\cite{kapovich}.

\subsection{Self-products of the Cartwright-Steger surface}\label{selfproductscs}

It is of course interesting to look for more examples of projective (or closed K\"ahler) manifolds endowed with an irrational pencil to which one can apply Theorems~\ref{bnestinfini} and~\ref{thgroup}. One way to build new examples is to use the {\it Cartwright-Steger surface}. Recall that this surface is a smooth compact complex surface which is a quotient of the unit ball $B$ of $\mathbb{C}^{2}$. We will denote it by $Y$ in what follows. Hence there exists a torsion-free cocompact lattice $\Gamma < {\rm PU}(2,1)$ such that $Y=B/\Gamma$. The surface $Y$ is characterized (up to changing the sign of the complex structure) by the fact that its Euler characteristic is equal to $3$ and its first Betti number is equal to $2$. It was discovered in~\cite{cs} in the context of the classification of fake projective planes. It was further studied by several authors, see for instance~\cite{cky0,stover,vidussi}. 

We now denote by $h : Y \to E$ the Albanese map of $Y$, whose target is an elliptic curve since $b_{1}(Y)=2$. Cartwright, Koziarz and Yeung~\cite{cky0} have proved that the map $h$ has isolated critical points and Koziarz and Yeung later proved that these critical points are nondegenerate~\cite{ky1}. We can thus consider the product $Y^{b}$ of $Y$ with itself $b$ times and the map
\begin{equation}\label{pencilsury}
h+\cdots + h : Y^b\to E.
\end{equation}
This provides natural examples to which one can apply Theorem~\ref{bnestinfini}. We discuss these examples in section~\ref{excs}. Denoting by $\Gamma < {\rm PU}(2,1)$ the fundamental group of the Cartwright-Steger surface, our construction together with Theorem~\ref{thgroup} immediately implies:

\begin{imain} The direct product of $b$ copies of $\Gamma$ contains a coabelian normal subgroup $N$ which is of type ${\rm FP}_{2b-1}$ but satisfies that $H_{2b}(N,\mathbb{Q})$ is not finitely generated. 
\end{imain}

The group $N$ appearing above is the kernel of the morphism on fundamental groups induced by the map~\eqref{pencilsury}. The fact that $N$ is of type ${\rm FP}_{2b-1}$ will be explained at the end of section~\ref{proof}, in Remark~\ref{fpness}. For a description of the lattice $\Gamma$, we refer the reader to~\cite{cky0}.


\section{Growth of the $n$-th Betti number}\label{proof} 

We prove here Theorem~\ref{bnestinfini}. We let $f : X \to S$ and $\widehat{X}\to X$ be as in section~\ref{ipabn}. Let $\widehat{S}$ be the universal cover of $S$ and $\widehat{f} : \widehat{X}\to \widehat{S}$ be a lift of $f$. The surface $\widehat{S}$ is topologically a plane. We now make the following:

\noindent {\bf Observation.} There exists a $C^{\infty}$ map $\widehat{f}_{0} : \widehat{X}\to \widehat{S}$, $C^{\infty}$-close to $\widehat{f}$, which is holomorphic in a neighborhood of its critical set and such that each critical point is nondegenerate. 

This observation is well-known, see~\cite{milnor}. To prove the existence of such a map $\widehat{f}_{0}$ we perturb $\widehat{f}$ in a neighborhood of each degenerate critical point. We first identify $\widehat{S}$ with $\C$ or the unit disc of $\C$. Let $q\in \widehat{X}$ be a critical point of $\widehat{f}$. We pick a neighborhood $U_q$ of $q$ that we identify with an open ball $B(0,\varepsilon)$ centered at the origin of $\C^{n}$ via the choice of some holomorphic coordinates. So we can consider the map $\widehat{f} : B(0,\varepsilon)\simeq U_q \to \widehat{S} \subset \C$. Let $\ell_q$ be a linear form on $\C^n$ which is a regular value of the map $d\widehat{f} : B(0,\varepsilon) \to (\C^{n})^{\ast}$ and which is small enough (if $q$ is nondegenerate, we take $\ell_q =0$). Then the map $\widehat{f}-\ell_q : B(0,\varepsilon)\to \widehat{S}\subset \C$ still takes values in $\widehat{S}$ and has finitely many critical points which are all nondegenerate. The number of its critical points is exactly the Milnor number of the critical point $q$ for the map $\widehat{f}$. If $\ell_q$ is small enough, we can assume that all critical points lie in the ball $B(0,\frac{\varepsilon}{2})$ (see Appendix B in~\cite{milnor}, in particular the remark on page 113, or~\cite[\S 5.4]{ebeling} for more details). In the following we write $V_q=B(0,\frac{\varepsilon}{2})$. We perform the previous construction for each point of the critical set $Crit(\widehat{f})$ of $\widehat{f}$, assuming that the open sets $$(U_q)_{q\in Crit(\widehat{f})}$$ are disjoint. One builds the map $\widehat{f}_{0}$ by declaring that $\widehat{f}_{0}$ is equal to $\widehat{f}-\ell_q$ on $V_q$, to $\widehat{f}$ outside the union of the open sets $(U_q)_{q\in Crit(\widehat{f})}$ and to a deformation between $\widehat{f}$ and $\widehat{f}-\ell_q$ on $U_q-V_q$. This can be performed in such a way that every critical point of $\widehat{f}_{0}$ is contained in one of the $V_q$'s. Although this is not necessary, we also observe that all this can be done in an equivariant way for the action of $\pi_{1}(X)$, so that the map $\widehat{f}_{0}$ descends to a map $f_{0} : X\to S$ homotopic to $f$. This implies the statement of the observation.

Since all the discussion below is topological, we will now work with the map $\widehat{f}_{0}$ instead of $\widehat{f}$. As in~\cite{dps}, we study the topology of $\widehat{X}$ by viewing it as the increasing union of a well-chosen sequence of compact subsets. We write 
$$\widehat{S}=\cup_{k\ge 1}D_{k}$$
where each $D_{k}$ is homeomorphic to a closed disk, $D_{k}\subset {\rm Int}(D_{k+1})$, no critical value of $\widehat{f}_0$ is contained in the boundary of $D_k$, $D_{k+1}-D_k$ contains exactly one critical value for $k\ge 1$ and $D_1$ contains no critical value. We set $\widehat{X}_{k}:=\widehat{f}_{0}^{-1}(D_k)$. Hence $\widehat{X}_1$ retracts onto a smooth fiber of $\widehat{f}_0$. For $k\ge 2$, the topology of the space $\widehat{X}_{k}$ is described thanks to the following proposition. 

\begin{prop}\label{vcht}
The space $\widehat{X}_{k+1}$ has the homotopy type of a space obtained from $\widehat{X}_{k}$ by gluing to it a finite number $m_k >0$ of $n$-dimensional cells. \
\end{prop}

This proposition is well-known (see e.g. Lemma 3.3 in~\cite{dps} for a related statement, although in that lemma the authors work with the original map $\widehat{f}$ instead of our perturbation $\widehat{f}_{0}$). We briefly sketch its proof below, for the reader's convenience. 

\noindent {\it Proof.} Let $c : [0,1] \to {\rm Int}(D_{k+1})$ be an embedded arc going from a boundary point $c(0)$ of $D_{k}$ to the unique critical value contained in $D_{k+1}-D_{k}$. We assume that $c(t)\notin D_{k}$ for $t>0$. Let $D^{\ast}$ be a small disk centered at $c(1)$ and contained in ${\rm Int}(D_{k+1} \setminus D_{k})$ (see Figure 1). 

\begin{figure}[h]
\centering
\includegraphics [scale=.65]{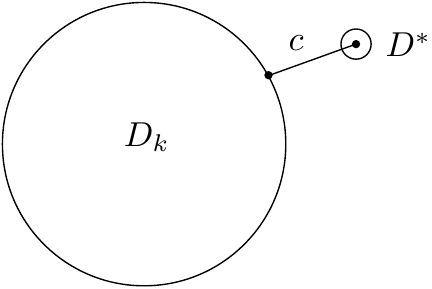}
\caption{The disks $D_k$ and $D^{\ast}$}\label{fig}
\end{figure}

Since the restriction of $\widehat{f}_0$ to the preimage of the set of regular values is a locally trivial fibration, $\widehat{X}_{k+1}$ deformation retracts onto $\widehat{X}_{k}\cup \widehat{f}_{0}^{-1}(c([0,1])\cup D^{\ast})$. Let $m_{k}>0$ be the number of critical points in the level $\widehat{f}_{0}^{-1}(c(1))$. Let $x_{1}, \ldots , x_{m_{k}}$ be the corresponding critical points. We fix a regular value $t=c(1-\delta)$ (for a small $\delta >0$) of $\widehat{f}_{0}$ close enough to $c(1)$. According to the theory of Lefschetz fibrations (see~\cite[\S 14.2.2]{voisin}) one can find $m_k$ disjoint $(n-1)$-dimensional spheres 
$$S_{1}\sqcup \ldots \sqcup S_{m_{k}}\subset \widehat{f}_{0}^{-1}(t)$$
(each sphere $S_j$ being contained in an arbitrary small neighborhood of $x_j$ fixed in advance) such that $\widehat{f}_{0}^{-1}(D^{\ast})$ has the homotopy type of the space obtained from $\widehat{f}_{0}^{-1}(t)$ by gluing an $n$-dimensional ball to each of the spheres $S_j$. By fixing a trivialization of the fibration
$$\widehat{f}_{0} : \widehat{f}_{0}^{-1}(c([0,1)))\to c([0,1)),$$ each sphere $S_j$ can be identified to a sphere $S_{j}^{\ast}\subset \widehat{f}_{0}^{-1}(c(0))$ in such a way that the $(S_{j}^{\ast})_{1\le j \le m_{k}}$ are disjoint. Hence $\widehat{X}_{k+1}$ retracts onto the space obtained from $\widehat{X}_{k}$ by gluing a ball to each of the spheres $S_{j}^{\ast}$. This proves the result.\hfill $\Box$

Proposition~\ref{vcht} immediately implies the following: 

\noindent {\bf Observation.} The sequence $b_{n-1}(\widehat{X}_{k})={\rm dim}_{\mathbb{Q}} H_{n-1}(\widehat{X}_{k},\mathbb{Q})$ is decreasing with $k$.    

Let $k_{0}$ be a large enough integer such that the sequence $$(b_{n-1}(\widehat{X}_{k}))_{k\ge k_{0}}$$ is constant. We will now prove that the sequence $(b_{n}(\widehat{X}_{k}))_{k\ge k_{0}}$ is strictly increasing with $k$ and that each map 
\begin{equation}\label{inclusho}
H_{n}(\widehat{X}_{k},\mathbb{Q})\to H_{n}(\widehat{X}_{k+1},\mathbb{Q})
\end{equation}
induced by the inclusion $\widehat{X}_{k}\hookrightarrow \widehat{X}_{k+1}$ is injective for $k\ge k_0$. This immediately implies Theorem~\ref{bnestinfini} since the group $H_{n}(\widehat{X},\mathbb{Q})$ is the direct limit of the $H_{n}(\widehat{X}_{k},\mathbb{Q})$'s.

 We use the notations from the proof of Proposition~\ref{vcht}. Let $k\ge k_{0}$. We know that $\widehat{X}_{k+1}$ has the homotopy type of a space $W_{k+1}$ obtained by gluing a ball to each of the spheres $$S_{1}^{\ast}\sqcup \ldots \sqcup S_{m_{k}}^{\ast}\subset \widehat{X}_{k}.$$ We write:
 \begin{equation}\label{pmv}
 W_{k+1}=\widehat{X}_{k}\cup B_{1} \cup \ldots \cup B_{m_{k}}
 \end{equation}
where each $B_{j}$ is homeomorphic to an $n$-dimensional ball, $B_{j}\cap B_l=\emptyset$ for $l\neq j$ and $\widehat{X}_{k}\cap B_{j}$ is equal to the boundary of $B_j$ (or to the sphere $S_{j}^{\ast}$ depending on whether one views it inside $\widehat{X}_{k}$ or $B_j$). Since the inclusion of $\widehat{X}_{k}$ into $\widehat{X}_{k+1}$ induces an isomorphism on $(n-1)$-dimensional homology groups, the same occurs for each inclusion $\widehat{X}_{k}\hookrightarrow W_{k+1}$. We now apply the Mayer-Vietoris exact sequence to the decomposition of $W_{k+1}$ given in~\eqref{pmv}.

We obtain (all homology groups being with $\mathbb{Q}$ coefficients): 

\begin{equation}
{\small \xymatrix{ & H_{n}(\sqcup_{j=1}^{m_{k}} \partial B_{j}) \ar[r] & H_{n}(\widehat{X}_{k})\oplus H_{n}(\sqcup_{j=1}^{m_{k}}B_{j}) \ar[r] & H_{n}(W_{k+1}) \ar[r] & \\
\ar[r] & H_{n-1}(\sqcup_{j=1}^{m_{k}} \partial B_{j}) \ar[r] & H_{n-1}(\widehat{X}_{k})\oplus H_{n-1}(\sqcup_{j=1}^{m_{k}}B_{j}) \ar[r] & H_{n-1}(W_{k+1}). &}}
\end{equation}
Since the groups $H_{n}(\sqcup_{j=1}^{m_{k}} \partial B_{j})$, $H_{n}(\sqcup_{j=1}^{m_{k}}B_{j})$ and $H_{n-1}(\sqcup_{j=1}^{m_{k}}B_{j})$ are zero we obtain:
\begin{equation}
{\small \xymatrix{\{0\} \ar[r] & H_{n}(\widehat{X}_{k}) \ar[r] & H_{n}(W_{k+1}) \ar[r] & H_{n-1}(\sqcup_{j=1}^{m_{k}} \partial B_{j}) \ar[r] & H_{n-1}(\widehat{X}_{k}) \ar[r] & H_{n-1}(W_{k+1}).}}
\end{equation}
The last arrow on the right being an isomorphism, this implies that the following sequence is exact:
\begin{equation}
{\small \xymatrix{\{0\} \ar[r] & H_{n}(\widehat{X}_{k}) \ar[r] & H_{n}(W_{k+1}) \ar[r] & H_{n-1}(\sqcup_{j=1}^{m_{k}} \partial B_{j}) \ar[r] & \{0\}.}}
\end{equation}
This implies that each inclusion $\widehat{X}_{k}\hookrightarrow W_{k+1}$ (and hence the inclusion $\widehat{X}_{k}\hookrightarrow \widehat{X}_{k+1}$) induces an injective map on $H_{n}(\cdot , \mathbb{Q})$ and that $b_{n}(\widehat{X}_{k})=b_{n}(\widehat{X}_{k+1})+m_{k}$. Note that $m_{k}>0$. This is the desired result and concludes the proof of Theorem~\ref{bnestinfini}.

\begin{rem}\label{fpness} Assume that $X$ is aspherical. Then the space $\widehat{X}$ is a $K(\Lambda,1)$. It has the homotopy type of a smooth fiber of $\widehat{f}$ (or $f$) with an infinite number of $n$-dimensional cells attached. This observation already appears in~\cite{dps} (see Corollary 5.8). This implies that $\Lambda$ is of type $\mathscr{F}_{n-1}$ (hence ${\rm FP}_{n-1}$), although it is not of type ${\rm FP}_{n}$.  
\end{rem}


\section{Examples built from the Cartwright-Steger surface}\label{excs}

We resume with the notation from Section~\ref{selfproductscs}. The Cartwright-Steger surface is denoted by $Y$ and $h : Y \to E$ is its Albanese map. Besides considering the products $Y\times \cdots \times Y$, one can also build more examples by combining the construction by Dimca, Papadima and Suciu and our construction. We fix a family of ramified covers $p_i : \Sigma_i \to E$ of the elliptic curve $E$ ($1\le i \le a$), where each $\Sigma_i$ has negative Euler characteristic. We then consider the map 
$$f : \Sigma_{1} \times \cdots \times \Sigma_{a} \times Y \times \cdots \times Y\to E$$
(where there are $b\ge 1$ copies of $Y$) which is the sum of the $p_i$'s and of the map $h$ on each copy of $Y$. All the results until the end of this section also apply when $a=0$, i.e. when one studies the map $f=h+\cdots +h$ as in~\eqref{pencilsury}.

The map $f$ has a finite non-empty set of critical points and connected fibers. This last point follows from the fact that $h : Y \to E$ has connected fibers. We denote by $\Lambda$ the kernel of the map induced by $f$ on fundamental groups. Theorem~\ref{dpsfirst} and Theorem~\ref{thgroup} imply that the group $H_{a+2b}(\Lambda,\mathbb{Q})$ is not finitely generated and that $\Lambda$ is projective if and only if $2b+a\ge 3$. The following proposition shows that the group $\Lambda$ is of a different nature compared to the examples from~\cite{dps,Llo-fourier}.

\begin{prop}\label{paspsg}
No finite index subgroup of $\Lambda$ embeds in a direct product of surface groups.
\end{prop}

By a {\it surface group} we mean here the fundamental group of an oriented surface of finite type (open or closed). Hence a surface group is either free or the fundamental group of a closed oriented surface. To prove Proposition~\ref{paspsg}, we will make use of the following theorem due to Bridson, Howie, Miller and Short~\cite{bhms}.

\begin{main}\label{tbmhs}
Let $F_1,...,F_m$ be surface groups.
Let $G$ be a subgroup of the direct product $F_1\times \cdots \times F_m$. If $G$ is of type ${\rm FP}_{m}$, then $G$ is virtually isomorphic to a direct product of the form $H_{1}\times \cdots \times H_{k}$ where $k\le m$ and each $H_i$ is a surface group. In particular $G$ is of type ${\rm FP}_{\infty}$. 
\end{main}

Besides this theorem, we will also use the fact that the property of being ${\rm FP}_m$ (for $m\in \mathbb{N}\cup \{\infty\}$) is invariant under passing to a subgroup or an overgroup of finite index~\cite[VIII.5.1]{brown}.

\noindent {\it Proof of Proposition~\ref{paspsg}.} The group $\Lambda$ sits inside the direct product 
\begin{equation}\label{dpp}
\pi_{1}(\Sigma_1)\times \cdots \times \pi_{1}(\Sigma_{a})\times \pi_{1}(Y)\times \cdots \times \pi_{1}(Y).
\end{equation}
The {\it factors} of this direct product are the subgroups of the form 
$$\{1\} \times \cdots \times \pi_{1}(\Sigma_{j})\times \cdots \times \{1\}$$
or
$$\{1\} \times \cdots \times \pi_{1}(Y)\times \cdots \times \{1\}.$$
Each factor intersects $\Lambda$ nontrivially. This implies that $\Lambda$ contains copies of the group $\mathbb{Z}^{a+b}$. The Gromov hyperbolicity of each factor implies that $\mathbb{Z}^{a+b+1}$ does not embed in $\Lambda$. Now suppose that a finite index subgroup $\Lambda_1$ of $\Lambda$ embeds in a direct product of surface groups $F_1\times \cdots \times F_m$. By taking $m$ to be minimal, we may assume that $L_i=\Lambda_{1}\cap F_i$ is nontrivial for $i=1,...,m$. Otherwise $\Lambda_1$ embeds in a direct product of $m-1$ surface groups. Since $\Lambda_1$ does not contain any nontrivial abelian normal subgroup, this implies that each $F_i$ is non-abelian, hence hyperbolic. A similar argument as before then shows that $\Lambda_1$ contains copies of $\mathbb{Z}^m$ but no copy of $\mathbb{Z}^{m+1}$. Hence $m=a+b$. By Remark~\ref{fpness}, the group $\Lambda$ is of type ${\rm FP}_{2b+a-1}$, hence $\Lambda_1$ is. Since $2b+a-1\ge m=a+b$, $\Lambda_{1}$ is of type ${\rm FP}_{m}$ and Theorem~\ref{tbmhs} implies that 
$\Lambda_{1}$ is of type ${\rm FP}_{\infty}$. This contradicts the fact that $\Lambda_1$ is not of type ${\rm FP}_{a+2b}$.\hfill $\Box$

Finally, we compute, in some cases, the first Betti number of the group $\Lambda=\ker(f_{\ast})$.  A similar computation appears in~\cite[\S 7]{Llo-17}, which applies to some of the examples built in~\cite{dps,Llo-fourier,Llo-17}.

\begin{prop}\label{fbn}
Assume that $a+2b\ge 3$. Assume furthermore that $b\ge 2$ or that $b=1$ and that the map $\pi_{1}(\Sigma_{j})\to \pi_{1}(E)$ is surjective for some $j\in \{1,\ldots , a\}$. Then the first Betti number of $\Lambda$ is equal to: $$b_1(\Sigma_1\times \cdots \Sigma_{a}\times Y^{b})-2.$$
\end{prop}

\noindent {\it Proof of Proposition~\ref{fbn}.} We consider the surjective homomorphism $$\Lambda \to \pi_{1}(\Sigma_1)\times \cdots \times \pi_{1}(\Sigma_{a})\times \pi_{1}(Y)^{b-1}$$ obtained by considering the inclusion of $\Lambda$ in the direct product~\eqref{dpp} and by projecting onto the first $a+b-1$ factors. Its kernel $N$ consists of elements of the form 
$$(1,\ldots , 1,g)\in \pi_{1}(\Sigma_1)\times \cdots \times \pi_{1}(\Sigma_{a})\times \pi_{1}(Y)\times \cdots \times \pi_{1}(Y)$$
where $g\in \ker(h_{\ast})$; it is isomorphic to $\ker(h_{\ast})$.  Hence we have the following exact sequence:
$$0\to N\to \Lambda \to \pi_{1}(\Sigma_1)\times \cdots \times \pi_{1}(\Sigma_{a})\times \pi_{1}(Y)^{b-1}\to 0.$$
It induces the following short exact sequence (see~\cite[VII.6]{brown}, all homology groups are taken with $\mathbb{Z}$ coefficients):
\begin{equation}\label{calculbetti}
H_{1}(N)_{\Lambda}\to H_{1}(\Lambda)\to H_{1}(\pi_{1}(\Sigma_1)\times \cdots \times \pi_{1}(\Sigma_{a})\times \pi_{1}(Y)^{b-1})\to 0.
\end{equation}
Here $H_{1}(N)_{\Lambda}$ is the group of coinvariants of $H_{1}(N)$ for the $\Lambda$-action. It is isomorphic to the quotient of $N$ by the group $[N,\Lambda]$ generated by commutators of elements of $\Lambda$ and of $N$. Note that if $x=(1, \ldots , 1, g)\in N$ and $y=(y_{1}, \ldots, y_{a}, h_{1},\ldots , h_{b})\in \Lambda$ then
$$xyx^{-1}y^{-1}=(1,\ldots , 1,gh_{b}g^{-1}h_{b}^{-1}).$$
Hence when we identify $N$ with $\ker(h_{\ast})$, $[N,\Lambda]$ is identified with the group $[\ker(h_{\ast}),\pi_{1}(Y)]$ (we are using here that $b\ge 2$ or that one of the $\pi_{1}(\Sigma_{j})$ surjects onto $\pi_{1}(E)$). In particular the groups $H_{1}(N)_{\Lambda}$ and $H_{1}(\ker(h_{\ast}))_{\pi_{1}(Y)}$ are isomorphic. 
Now the short exact sequence
$$0\to \ker (h_{\ast}) \to \pi_{1}(Y) \to \mathbb{Z}^{2}\to 0$$
induces the short exact sequence (see~\cite[VII.6]{brown} again):
$$H_{2}(\mathbb{Z}^{2})\to H_{1}(\ker(h_{\ast}))_{\pi_{1}(Y)}\to H_{1}(Y)\to H_{1}(\mathbb{Z}^{2})\to 0.$$
Since the map $H_{1}(Y)\otimes \mathbb{Q} \to H_{1}(\mathbb{Z}^{2})\otimes \mathbb{Q}$ is an isomorphism, we obtain that $H_{1}(\ker(h_{\ast}))_{\pi_{1}(Y)}\otimes \mathbb{Q}$ has dimension at most $1$. Hence $H_{1}(N)_{\Lambda}\otimes \mathbb{Q}$ has dimension at most $1$. Since $\Lambda$ is K\"ahler and hence has even first Betti number, this implies that the first arrow in~\eqref{calculbetti} has finite image. Hence $H_{1}(\Lambda)\otimes \mathbb{Q}$ and $H_{1}(\pi_{1}(\Sigma_{1})\times  \cdots \times \pi_{1}(\Sigma_{a})\times \pi_{1}(Y)^{b-1})\otimes \mathbb{Q}$ are isomorphic. This gives the desired result.\hfill $\Box$

\begin{rem} The topology of the generic fiber of our map $f$ (mainly its homotopy groups up to dimension $a+2b-1$) can be described thanks to Theorem 5.2 and Corollary 5.4 from~\cite{dps}. Indeed all the hypotheses required in loc. cit. are met in our context, except possibly for the hypothesis on resonance varieties. However, that hypothesis is only used in~\cite{dps} to prove that the group $H_{a+2b}(\Lambda,\mathbb{Q})$ is not finitely generated, a result that follows from our Theorem~\ref{thgroup}. 
\end{rem}

\begin{rem} By considering the case where $a\in \{0,1\}$, our construction provides for each $n\ge 2$ an example of a ${\rm CAT}(0)$ group $G$ containing a subgroup of type ${\rm FP}_{n-1}$ but not ${\rm FP}_{n}$ and such that $G$ does not contain free Abelian subgroups of rank greater than $\left \lfloor \frac{n+1}{2}\right\rfloor$. See~\cite{brady,rkropholler} for related results and motivation. The article~\cite{rkropholler} produces other examples with a smaller bound on the rank of Abelian subgroups. More precisely, for each positive integer $n$, Kropholler gives in~\cite{rkropholler} an example of a group of type $\mathscr{F}_{n-1}$ but not of type $\mathscr{F}_{n}$, which does not contain free Abelian groups of rank greater than $\left \lceil \frac{n}{3}\right \rceil$ and which is a subgroup of a ${\rm CAT}(0)$ group. 
\end{rem}

{\bf Acknowledgements.} We would like to thank Claudio Llosa Isenrich and Valdo Tatitscheff for their comments on this text. We would also like to thank the referee for their constructive comments on the text.

\bigskip
\bigskip
\begin{small}
\begin{tabular}{llll}
Francisco Nicol\'as & & & Pierre Py\\
IRMA & & & IRMA \\
Universit\'e de Strasbourg \& CNRS & & & Universit\'e de Strasbourg \& CNRS \\
67084 Strasbourg, France & & & 67084 Strasbourg, France \\
fnicolascardona@math.unistra.fr & & & ppy@math.unistra.fr\\    
\end{tabular}
\end{small}

\end{document}